\documentclass[final,numberedheadings]{aipproc}

\layoutstyle{6x9}

\usepackage{amsmath,amssymb,amsfonts,amsthm,amscd}
\usepackage[mathscr]{eucal}

\newtheorem{thm}{Theorem}
\newtheorem{prop}{Proposition}

\newtheorem*{coro}{Corollary}
\theoremstyle{definition}
\newtheorem{de}{Definition}%[section]
\newtheorem*{Rem}{Remark}

\newcommand{\co}{\colon\thinspace}

\renewcommand{\geq}{\geqslant}

\DeclareMathOperator{\ad}{ad}

\DeclareMathOperator{\Vect}{\mathfrak{X}}
\DeclareMathOperator{\Vectn}{\Vect_{-\vphantom{1}}}
\DeclareMathOperator{\Vectp}{\Vect_{+}}

\newcommand{\w}{{\mathbf{w}}}

\newcommand{\der}[2]{{\frac{\partial {#1}}{\partial {#2}}}}
\newcommand{\lder}[2]{{\partial {#1}/\partial {#2}}}

\newcommand{\Z}{{\mathbb Z_{2}}}
\newcommand{\ZZ}{{\mathbb Z}}

\newcommand{\p}{\partial}
\newcommand{\fun}{C^{\infty}}

\newcommand{\la}{{\lambda}}
\newcommand{\x}{{\xi}}

\newcommand{\ft}{{\tilde f}}
\newcommand{\jt}{{\tilde j}}

\newcommand{\gt}{{\tilde g}}

\newcommand{\jtt}{{\tilde j}}

\newcommand{\itt}{{\tilde i}}

\newcommand{\ut}{{\tilde u}}
\newcommand{\vt}{{\tilde v}}

\newcommand{\mut}{{\tilde \mu}}

\begin{document}

\title{$Q$-manifolds and Higher Analogs of Lie Algebroids}


\author{Theodore~Th.~Voronov}{ address={School of Mathematics, The
    University of Manchester, Oxford Road, Manchester, M13 9PL, UK} }

\begin{abstract}
We show how the relation between $Q$-manifolds and  Lie algebroids   extends to   ``higher'' or ``non-linear'' analogs of Lie algebroids. We study the identities satisfied by  a new algebraic structure that arises as a replacement of operations on sections of a  Lie algebroid. When the base is a point, we obtain  a generalization of Lie superalgebras.   %We also consider  a particular case of transitive  ``non-linear''  Lie algebroids and show an application  of  the non-Abelian Poincar\'{e} lemma introduced earlier by the author.
\end{abstract}

\maketitle

%%%%%%%%%%%%%%%%%%%%%%%%%%%%%%%%%%%%%%%%%%%%
%% MAINMATTER
%%%%%%%%%%%%%%%%%%%%%%%%%%%%%%%%%%%%%%%%%%%%

\section{Recollection: Lie algebroids in the language of $Q$-manifolds}

In the sequel, $\Vect(M)$ denotes the space of all vector fields on a (super)manifold $M$. In the graded case, see below, we shall use the same notation for the subspace of graded vector fields, i.e., the direct sum of the homogeneous vector fields.

Let us recall the following definition.

\begin{de} A smooth vector bundle $E\to M$ is a \emph{Lie algebroid} if it is endowed with the following additional structure: a Lie bracket on sections and a vector bundle map $a\co E\to TM$ (called the \emph{anchor}) so that the following identity is satisfied:
\begin{equation*}
    [u,fv]=a(u)f\, v +(-1)^{\ut\ft}f[u,v]\,,
\end{equation*}
for arbitrary smooth sections of $E\to M$ and functions on $M$.
\end{de}

Comment: in this definition both $E$ and $M$ are  supermanifolds, so we have actually defined a super Lie algebroid. It makes no principal difference with the Lie algebroids defined in the setting of ordinary manifolds. The standard source about Lie algebroids (and Lie groupoids) is  the fundamental book by Mackenzie~\cite{mackenzie:book2005}.

Since $E\to M$ is a vector bundle, we may endow it with a natural $\ZZ$-graded structure where the functions on the base have weight $0$ and linear functions on the fibers have weight $1$. It induces the same  grading on the reversed parity bundle $\Pi E$, so if we denote local coordinates on $M$ by $x^a$ and the linear coordinates on the fibers of $\Pi E$ corresponding to a local frame $e_i$ for $E$ by $\x^i$, then $\w(x^a)=0$, $\w(\x^i)=1$\,. (Here $\w$ denotes weight.) We should warn the reader that in general there is no link between parity and weight. See~\cite{tv:graded} for more details on graded manifolds.

Recall that a \emph{$Q$-manifold} or a manifold with a \emph{$Q$-structure} is a supermanifold  endowed with a homological vector field, i.e., an odd vector field $Q$ such that $[Q,Q]=2Q^2=0$\,. There may be  an extra $\ZZ$-grading and the field $Q$ may be homogeneous w.r.t. weight. (In local coordinates, the condition $Q^2=0$ reads $Q^b\p_bQ^q=0$ if $Q=Q^a(x)\p_a$\,.)

\begin{prop}[Vaintrob~\cite{vaintrob:algebroids}] A Lie algebroid structure on $E$ is equivalent to a $Q$-structure  on $\Pi E$ of weight $+1$\,:
\begin{equation*}
    Q=\x^iQ_i^a(x)\,\der{}{x^a}+\frac{1}{2}\,\x^i\x^j Q_{ji}^k
    (x)\,\der{}{\x^k}\,.
\end{equation*}
\end{prop}

For a given $Q\in \Vect (\Pi E)$ as above,  the   bracket of sections of $E$ and the anchor $a\co E\to TM$ are  defined by the formulas
\begin{equation*}
    a(e_i)=Q_i^a(x)\,\der{}{x^a}\,,
\end{equation*}
and
\begin{equation*}
    [e_i,e_j]=(-1)^{\jtt} Q_{ij}^k(x)\,e_k\,
\end{equation*}
(on basis sections). The equation $Q^2=0$ contains all the axioms of a Lie algebroid in a compact form. This can be checked directly.

The relation between the operations on sections of $E$ and the vector field $Q$ on $\Pi E$ can be defined in a more invariant fashion as follows. The space of vector fields on $\Pi E$ is $\ZZ$-graded by weights, so that
\begin{equation*}
    \Vect (\Pi E)=\bigoplus_{n\geq -1} \Vect_n(\Pi E)\,.
\end{equation*}
There are no nonzero vector fields of weights less than $-1$, and the vector fields of weight $-1$ are in a one-to-one correspondence with the sections of $E$. We define an odd linear map
\begin{equation*}
    i\co \fun(M,E)\to \Vect(\Pi E)\,, \quad u=u^ke_k\mapsto i_u=(-1)^{\ut}u^k\der{}{\x^k}\,,
\end{equation*}
which   is an isomorphism onto $\Vect_{-1}(\Pi E)$. Now the relation between the homological vector field $Q\in \Vect_{1}(\Pi E)$ and the Lie algebroid structure on $E$ is given by the derived bracket formulas
\begin{equation*}
    a(u)f=\bigl[[Q,i_u)],f\bigr]\,,
\end{equation*}
\begin{equation*}
    i_{[u,v]}=(-1)^{\ut}\bigl[[Q,i_u],i_v\bigr]\,.
\end{equation*}
See, for example, ~\cite{tv:graded}. This is the derived bracket of $Q$ restricted on  the subalgebra $\fun(M)\oplus i(\fun(M,E))$ in the Lie algebra $\fun(M)\oplus \Vect(\Pi E)$\,\footnote{Here the Lie bracket between functions is set to zero and, between a function and a vector field, it is the action of vector fields on functions.}. The required properties of the bracket of sections  and the anchor follow from the general properties of   derived brackets (see~\cite{yvette:derived, yvette:derived2}, and \cite{tv:graded}).

As it is known, one of the advantages of describing Lie algebroids in the language of $Q$-manifolds is the simplification of the description  of morphisms~\cite{vaintrob:algebroids}. We recall that briefly. A \emph{$Q$-manifold morphism} between $Q$-manifolds $(M_1,Q_1)$ and $(M_2,Q_2)$ is a smooth map $f\co M_1\to M_2$ such that the vector fields $Q_1$ and $Q_2$ are $f$-related. In local coordinates, if we denote coordinates on $M_1$ by $x^a$ and on $M_2$ by $y^i$, and the map $f$ is given by $y^i=y^i(x)$, then the condition that it is a  {$Q$-manifold morphism} is expressed by the equation
\begin{equation*}
    Q_2^i(y(x))=Q_1^a(x)\,\der{y^i}{x^a} \,,
\end{equation*}
in a self-explaining notation. Such maps can be obviously combined, so they make a category. Now, for Lie algebroids $E_1\to M_1$ and $E_2\to M_2$ over different bases,  a morphism of vector bundles $E_1\to E_2$ is a Lie algebroid morphism (as defined in~\cite{mackenzie:book2005}; a non-obvious notion) if and only if the induced map $\Pi E_1\to \Pi E_2$ is a $Q$-manifold morphism, which is much easier to handle.

This language is also very helpful for studying Lie bialgebroids and their `doubles', see~\cite{mackenzie:drinfeld}, \cite{roytenberg:thesis}, \cite{tv:graded}. Another advantage is the simplicity of introduction of  arbitrary double and multiple Lie algebroids~\cite{tv:mack}.

\section{Graded $Q$-manifolds and ``higher'' Lie algebroids. A new algebraic structure}

The philosophical idea is that a  $\ZZ$-grading should be viewed as a replacement of a linear structure when the latter is not available. Notice that for an algebra of polynomials, if all its generators are assigned the same weight (i.e., if weights are counted  as the usual degrees), then the automorphisms preserving weights are just the linear transformations of variables, so this is equivalent to a choice of   linear structure. If, however, different generators are assigned different weights, then among weight-preserving automorphisms there arise non-linear transformations.

Such a situation materializes, for example, for multiple vector bundles, say, double vector bundles  such as $TE$ or $T^*E$ for a given vector bundle $E\to M$. They are not vector bundles over the original base $M$ (no linear structure) and considering `total weight' for them leads to coordinates having weights $0, 1$, and $2$.

Along these lines we obtain a ``non-linear'' generalization of Lie algebras and Lie algebroids.
We need to recall the structure of graded manifolds given in~\cite{tv:graded}.

Let $F$ be a graded manifold. Denote local coordinates by $x^a$, $y^i$ where $\w(x^a)=0$ and $\w(y^i)\neq 0$.  There is an underlying non-graded supermanifold $F_0$ %(all functions on $F_0$ have zero weight)
embedded into $F$ as a closed submanifold and specified locally by the equations $y^i=0$.  As local coordinates on $F_0$ one may take the restrictions of $x^a$. Suppose now that $F$ is non-negatively graded, i.e., there are no coordinates of negative weights.  Then there is also a canonical projection $F\to F_0$, so that the graded manifold $F$ turns out to be  a fiber bundle over $M=F_0$. More precisely, there is  a tower of fibrations
\begin{equation}\label{eq.tower}
    F=F_{N}\to F_{N-1}\to \ldots \to F_{2}\to F_{1}\to F_0=M\,
\end{equation}
obtained by ordering  the local coordinates by increasing weights; the number $N$ is the maximal weight of local coordinates. The first fibration $F_1\to M$ is a vector bundle, all the higher fibrations $F_{k+1}\to F_k$ are `affine bundles' (i.e., the transition functions are affine transformations where the  linear homogeneous parts are the transformations between the coordinates of weight $k+1$ and the additive terms are polynomial expressions of total weight $k+1$ in the coordinates of smaller weights). So for the fiber bundle $F\to M$, the transition functions  are particular polynomial functions of fiber coordinates.  For $N=1$ we come back to vector bundles, and for $N>1$,  non-negatively graded manifolds should be viewed as the non-linear generalizations of such. (The number $N$ is sometimes referred to as the \emph{degree} of a non-negatively graded manifold.)

Suppose there is a homological vector field $Q$ on an arbitrary graded manifold $F$. As a section $F\to \Pi TF$, the odd vector field $Q$ itself can be tautologically regarded as analogous to the anchor of Lie a algebroid. Moreover, one can see that $Q^2=0$ is equivalent to the map $Q\co F\to \Pi TF$ being a $Q$-manifold morphism between $(F,Q)$ and $(\Pi TF,d)$\,. (Indeed, in local coordinates the map $Q$ is given by $(z^{\mu})\mapsto\bigl(z^{\mu}=z^{\mu}, dz^{\mu}=Q^{\mu}(z)\bigr)$, here $z^{\mu}$ denote the whole bunch of local coordinates on $F$,
so the condition that it is a $Q$-morphism, see the previous section,  amounts to the two equations: $dz^{\mu}=Q^{\mu}(z)$, which reproduces the definition of the map, and $0=Q^{\nu}\p_{\nu}Q^{\mu}$, which is exactly $Q^2=0$.) If $F$ is non-negatively graded, so we have a fiber bundle $p\co F\to M$, a closer analog of Lie algebroid anchor is the map $a=a_Q\co F\to \Pi TM$ defined as the composition $a_Q:=\Pi Tp\circ Q$\,. If in the local coordinates $x^a$, $y^i$, where   $\w(x^a)=0$ and $\w(y^i)> 0$,
\begin{equation*}
    Q=Q^a(x,y)\der{}{x^a}+Q^i(x,y)\der{}{y^i}\,,
\end{equation*}
then
\begin{equation*}
    a_Q\co (x^a,y^i)\mapsto \bigl(x^a,dx^a=Q^a(x,y)\bigr)\,.
\end{equation*}
Now, if we require that $a_Q\co F\to \Pi TM$ preserves weight and the vector field $Q$ is homogeneous, then $\w(Q)=1$ by necessity, in a close analogy with the Lie algebroid case. We have therefore arrived at the following definition.
\begin{de}[\cite{tv:graded}] \label{def.nonlinalg} A \emph{non-linear} or \emph{higher Lie algebroid} with base $M$ is a non-negatively graded manifold $F\to M$ where $M=F_0$ is endowed with a homological vector field $Q\in \Vect(F)$ of weight $+1$\,. The \emph{anchor} is the map $a=a_Q\co F\to \Pi TM$ defined above.
\end{de}
(For a Lie algebroid $E$, we have $F=\Pi E$,  but in the genuine non-linear case we cannot reverse parity, so there is no  ``$E$'' for $F$\,. In~\cite{tv:graded} a slightly more general setting was considered as ``non-linear Lie bialgebroids'':  a graded manifold without the non-negativity assumption endowed with a homological vector field of arbitrary weight $q$ and a compatible Poisson or Schouten bracket of weight $s$.)

It is helpful to understand a bit better the form of the field $Q\in\Vect(F)$ for a non-linear Lie algebroid $F$. The condition $\w(Q)=1$ implies that
\begin{equation}\label{eq.qfield}
    Q= y^iQ^a_i(x)\der{}{x^a}+Q^k(x,y)\der{}{y^k}\,, %\quad \text{where} \quad \w(Q^k)=\w(y^k)+1\,.
\end{equation}
where in the first term the coefficients are non-zero only for $\w(y^i)=1$ and in the second term $\w(Q^k)=\w(y^k)+1$.
Therefore the bundle map $a\co F\to \Pi TM$ is given, in local coordinates, by the formula
\begin{equation}
    dx^a=\sum_{\w(y^i)=1}y^iQ^a_i(x)\,.
\end{equation}
We see that it factors through a vector bundle morphism $a'\co F_1\to \Pi TM$, where the vector bundle $F_1\to M$ is the first fibration in the tower~\eqref{eq.tower}.

Now we shall describe an algebraic structure corresponding to a non-linear Lie algebroid. Unlike the classical case, we cannot use sections of the bundle $F\to M$, because it is not a vector bundle and its sections cannot be added or multiplied by functions. However we can consider instead graded vector fields on the total space $F$. Suppose the maximal weight of local coordinates on $F$ equals $N$. We have
\begin{equation*}
    \Vect (F)=\bigoplus_{n\geq -N} \Vect_n(F)\,.
\end{equation*}
The graded Lie superalgebra $\Vect (F)$ can be decomposed into the sum of two subalgebras:
\begin{equation*}
    \Vect (F)=\Vectn(F)\oplus \Vectp(F) \quad \text{where} \quad \Vectn(F)=\bigoplus\limits_{n= -N}^{n=-1}\Vect_n(F)\,,
    \Vectp(F)=\bigoplus\limits_{n\geq 0}\Vect_n(F)\,.
\end{equation*}
For ordinary Lie algebroids,  we have $\Vectn(\Pi E)=\Vect_{-1}(\Pi E)$; it is an  Abelian  subalgebra  isomorphic to the space of sections $\fun(M,E)$. Here  the subalgebra  $\Vectn(F)$ is   nilpotent; we shall use it as a replacement of the space of sections of a Lie algebroid. %We use it as a replacement of the space $\fun(M,E)$.
Let $P$ denote the projector onto $\Vectn(F)$. We  define new operations on the space $\Vectn(F)$ by the higher derived bracket formula~\cite{tv:higherder}\,:
\begin{equation*}
    \{u_1,\ldots,u_k\}_Q:=P[\ldots[[Q,u_1],u_2],\ldots,u_k]\,,
\end{equation*}
for $u_i\in \Vectn(F)$\,. Here $k=0,1,2,\ldots\,$. Since $\w(Q)=1$ and $\w(u_i)<0$ for all $u_i$\,, we immediately see the following. The operations are automatically odd and of weight $+1$. There is no $0$-ary bracket: $\{\varnothing\}_Q=P(Q)=0$; the unary bracket
 \begin{equation}\label{eq.defd1}
    d_Qu:=\{u\}_Q=P[Q,u]
 \end{equation}
is zero on $\Vect_{-1}(F)$ and equals $Du:=[Q,u]$ on $\Vect_{n}(F)$ for $n<-1$. The binary bracket does not really require the projector:
\begin{equation}\label{eq.defbracket1}
    \{u,v\}_Q=[[Q,u],v]\,,
\end{equation}
because the r.h.s. is automatically in $\Vectn(F)$ for $u,v\in \Vectn(F)$. Therefore all  the  higher brackets are compositions of the binary bracket $\{u,v\}_Q$ with the ``old'' brackets (the commutators of vector fields):
\begin{equation*}
    \{u_1,u_2,u_3,\ldots,u_k\}_Q=[\ldots [\{u_1,u_2\}_Q,u_3],\ldots,u_k]\,.
\end{equation*}
So it all boils down to studying the space $\Vectn(F)$ endowed with the three operations: the original commutator $[u,v]$, the differential $du$, and the new bracket $\{u,v\}$ (where for brevity we have suppressed the superscripts $Q$).

The properties of these operations can be summarized as follows. Since they can be obtained simply by counting weights and there is nothing peculiar for $\Vect(F)$, they hold in a more general setup. Consider   an arbitrary $\ZZ$-graded Lie superalgebra\,\footnote{We recall that the $\ZZ$-grading (weight) and the $\Z$-grading (parity) are, in general, not related.} $L=\bigoplus L_n$ and   define its decomposition as
\begin{equation*}
    L=V\oplus K \quad \text{where} \quad V=\bigoplus_{n<0}L_n\,, \ K=\bigoplus_{n\geq 0}L_n\,.
\end{equation*}
Let $Q$ be an odd element in $L_1$ and consider the corresponding inner derivation $D:=\ad Q$. Alternatively, consider an arbitrary odd derivation $D$ of weight $1$ on $L$. Define the derived brackets on $V$ by the same formula as above:
\begin{equation*}
    \{u_1,\ldots,u_k\}:=P[\ldots[Du_1,u_2],\ldots,u_k]\,,
\end{equation*}
where $P$ is the projector onto $V$ with the kernel $K$. Again, everything boils down to the operation $d$,
\begin{equation}\label{eq.defd}
    du = P(Du)\,,
\end{equation}
so    $du=Du$ if $\w(u)<-1$ and $du=0$ if $\w(u)=-1$,  and the binary bracket $\{\_,\_\}$,
\begin{equation} \label{eq.defbracket}
    \{u,v\}=[Du,v]\,.
\end{equation}
%Both new operations on $V$ are odd and of weight $1$.
Rename $L_n=:V_n$ for $n<0$ so that $V=\bigoplus_{n<0} V_n$.  We arrive at the following picture.
\begin{thm} The space $V$ is a negatively graded Lie superalgebra w.r.t. the original bracket $[\_,\_]$ (which is, in particular, even and of weight $0$), and it is nilpotent if the grading of $L$ is bounded from below.
The operations $d$ and $\{\_,\_\}$ defined by~\eqref{eq.defd} and~\eqref{eq.defbracket} respectively are both odd and of weight $+1$. We have
\begin{equation}\label{eq.complex}
    \ldots  \xrightarrow{\ d\ }V_{-n-1}\xrightarrow{\ d\ } V_{-n} \xrightarrow{\ d\ } \ldots \xrightarrow{\ d\ } V_{-1} \xrightarrow{\ d\ } 0\,.
\end{equation}
If $D^2=0$, then~\eqref{eq.complex} is a complex:
\begin{equation}\label{eq.dsq}
    d^2=0\,,
\end{equation}
and the following identities are satisfied:
\begin{align}
    d[u,v]&=\{u,v\}-(-1)^{\ut\vt}\{v,u\}\,, \label{eq.dbrack}\\
    d\{u,v\}&=-(-1)^{(\ut+1)\vt}\{v,du\}+(-1)^{\ut+1}\{u,dv\}+(-1)^{\ut}[du,dv]\,, \label{eq.dnewbrack}\\
    \{u,[v,w]\}&=[\{u,v\},w]+(-1)^{(\ut+1)\vt}[v,\{u,w\}]\,, \label{eq.leib} \\
    \{u,\{v,w\}\}&=(-1)^{\ut+1}\{\{u,v\},w\}+(-1)^{(\ut+1)(\vt+1)}\{v,\{u,w\}\}\,. \label{eq.jac}
\end{align}
\end{thm}
\begin{coro} For a non-linear Lie algebroid $F\to M$  there is a complex
\begin{equation}\label{eq.complex2}
   0\xrightarrow{\ d\ }  \Vect_{-N}(F)\xrightarrow{\ d\ }\Vect_{-N+1}(F)  \ldots \xrightarrow{\ d\ } \Vect_{-1}(F) \xrightarrow{\ d\ } 0\,.
\end{equation}
The original commutator of vector fields $[\_,\_]$, and the operations $d=d_Q$ and  $\{\_,\_\}=\{\_,\_\}_Q$ defined on $\Vectn(F)$  by~\eqref{eq.defd1} and~\eqref{eq.defbracket1} satisfy  the identities~\eqref{eq.dbrack}, \eqref{eq.dnewbrack}, \eqref{eq.leib}, and \eqref{eq.jac}\,.
\end{coro}

A few comments. Equation~\eqref{eq.jac} is the usual Jacobi identity (in the Leibniz form) for a derived bracket (see~\cite{yvette:derived}; we use a different sign convention, see~\cite{tv:graded, tv:higherder}). Equation~\eqref{eq.dbrack} is also familiar from~\cite{yvette:derived}. Note however that the odd bracket $\{\_,\_\}$ is \emph{not} derived from $d$; it is not possible to recover $\{\_,\_\}$ from $[\_,\_]$ and $d$ on $V$ (without the knowledge of the whole $L=V\oplus K$ and the derivation $D$ on $L$). Nevertheless, on $[V,V]$ the operators $d$ and $D$ coincide, and that is why~\eqref{eq.dbrack} holds. The differential $d$ is not a derivation of $[\_,\_]$ or $\{\_,\_\}$. The identity~\eqref{eq.dbrack} replaces the derivation property w.r.t. $[\_,\_]$ and at the same time measures the deviation of $\{\_,\_\}$ from (anti)symmetry. The identity~\eqref{eq.dnewbrack} is a `modified' form of the derivation property taking into account the lack of symmetry for $\{\_,\_\}$. The identity~\eqref{eq.leib} is the single ``mixed'' Jacobi--Leibniz identity involving both brackets  (an attempt to obtain a Leibniz-type rule for $[u,\_]$ w.r.t. $\{\_,\_\}$ leads to the same identity).  Notice also that the algebraic structure described by~\eqref{eq.dbrack}, \eqref{eq.dnewbrack}, \eqref{eq.leib}, and \eqref{eq.jac} degenerates into a Lie superalgebra structure for $N=1$ (more precisely, $[\_,\_]$ and $d$ disappear, and $\{\_,\_\}$ becomes a Lie bracket w.r.t. reversed parity after an inessential change of   sign).

So far we have not used the bundle structure of $F$ over $M$ or the anchor. The space $\Vectn(F)$ is not closed under the multiplication by functions on $F$ (because it can increase weight), but it is a graded module over $\fun(M)$. If we consider everything locally, it is clear that we obtain a locally-free graded module of finite rank over the sheaf of functions on $M$, so $\Vectn(F)$ can be identified with the space of sections of a certain graded vector bundle over $M$\,! (Note for comparison that $\fun(F)$ as well as the whole space $\Vect(F)$ are infinite-dimensional over $\fun(M)$.) We shall explore it further elsewhere, but now note the following. Since for $u\in \Vectn(F)$ and $f\in \fun(M)$, we have $uf=0$ by counting weights (so the elements of $\Vectn(F)$ are vertical), and so we have
\begin{equation}\label{eq.linearlie}
    [fu,v]=f[u,v]
\end{equation}
for all $f\in \fun(M)$ and $u,v\in \Vectn(F)$\,. Therefore we have a bundle of Lie (super)algebras, which is moreover a Lie superalgebra bundle (one has to examine the brackets of basis elements to see that). It remains to consider the relation of the module structure over $\fun(M)$ and the bracket $\{\_,\_\}$.

Define the \emph{anchor on the elements of $\Vectn(F)$} (this is in a new sense, compared to the map $a\co F\to \Pi TM$ considered above) as an odd linear transformation $\Vectn(F) \to \Vect(M)$ of weight $+1$ sending $u\in \Vectn(F)$ to $a(u)\in \Vect(M)$ given by
\begin{equation}\label{eq.anchor2}
    a(u)f:= [[Q,u],f]\,,
\end{equation}
where $f\in \fun(M)$. Then, from the definitions, we have the Leibniz identity
\begin{equation}\label{eq.leibanchor}
    \{u,fv\}=a(u)f\cdot v +(-1)^{(\ut+1)\ft}f\{u,v\}\,,
\end{equation}
for all $f\in \fun(M)$ and $u,v\in \Vectn(F)$\,. What is the relation of this ``algebraic'' anchor with the ``geometric'' anchor considered before?
\begin{prop} The two notions of anchor\,---\,defined as the bundle map $a\co F\to \Pi TM$, $a=\Pi Tp\circ Q$, and as the linear transformation $a\co \Vectn(F) \to \Vect(M)$ given by~\eqref{eq.anchor2}\,---\,are equivalent.
\end{prop}

To see that,  observe first that the transformation $u\mapsto a(u)$ is linear over $\fun(M)$, i.e.,
\begin{equation}\label{eq.anchorlin}
    a(gu)=(-1)^{\gt}ga(u)\,,
\end{equation}
for all $u\in \Vectn(F)$ and $g\in \fun(M)$\,. This follows by counting weights, but can also be seen from an explicit formula:
\begin{equation}
    a(u)=(-1)^{\ut+1} \sum_{\w(y^i)=1}u^iQ_i^a\der{}{x^a} \quad \text{where} \quad u=u^i\der{}{y^i}\,.
\end{equation}
If we denote a graded vector bundle over $M$ corresponding to $\Vectn(F)$ by $E$, we have an odd vector bundle morphism $E\to TM$. Also, because it is of weight $1$ and all elements of $TM$ have weight $0$, it vanishes on all $E_{n}$ except for $n=-1$ and, even more, it vanishes (locally) on all products $p^i(y)\lder{}{y^i}\in E_{-1}$ if $\w(y^i)\neq 1$. Hence it factors through the quotient by these products; the quotient vector bundle has the partial derivatives $\lder{}{y^i}|_{y=0}$ where $\w(y^i)=1$ as a local basis. It can be identified with the vector bundle $F_1$. So the  anchors defined in the two different ways can be identified  as  linear maps $F_1\to \Pi TM$.
In local coordinates,  the anchor in both senses is specified by the coefficients $Q^a_i(x)$ in the first term in~\eqref{eq.qfield}, exactly as for ordinary Lie algebroids.

\begin{Rem} Since the bracket $\{\_,\_\}$ is not symmetric, it is legitimate to ask whether a derivation property w.r.t. multiplication by functions on $M$ holds for the first argument. It is not hard to see that the identity is more complicated. Namely, for all $u,v\in \Vectn(F)$ and all $f\in \fun(M)$,
\begin{equation}\label{eq.leibanchor2}
    \{fu,v\}=(-1)^{\ft}f\{u,v\}+(-1)^{(\ft+\ut)\vt}a(v)f\cdot u +\delta f \bigl([u,v]\bigr)\,,
\end{equation}
where we have denoted by $\delta f$ the commutator of $d$ and the multiplication by $f$\,:
\begin{equation}\label{eq.deltaf}
    \delta f=[d,f]\,.
\end{equation}
It is an odd operator of weight $+1$ on $\Vectn(F)$, so it is zero on $\Vect_{-1}(F)$ and on $\Vect_{n}(F)$, $n<-1$, it can be expressed in coordinates as
\begin{equation*}
    \delta f =Qf= \sum_{\w(y^i)=1}y^iQ_i^a(x)\der{f}{x^a}
\end{equation*}
(a multiplication operator).
\end{Rem}

The identities~\eqref{eq.dsq}, \eqref{eq.dbrack}, \eqref{eq.dnewbrack}, \eqref{eq.leib}, and \eqref{eq.jac}, together with \eqref{eq.linearlie}, \eqref{eq.anchorlin}, and \eqref{eq.leibanchor}, completely describe the algebraic structure generated on $\Vectn(F)$ regarded as a graded module over $\fun(M)$\,. It is not difficult to axiomatize this and obtain a definition of what may be called a `two-layer Lie algebra'  (if we ignore the module structure and the anchor) or a `two-layer Lie algebroid / two-layer Lie pseudoalgebra'. By the ``two layers'' we mean the following. On the first layer there is a negatively-graded Lie superalgebra (over some ground commutative superalgebra). On the second layer there arise an odd differential $d$, an odd bracket $\{\_,\_\}$, and an odd anchor (transformation into derivations of the ground algebra), all of weight $+1$,  satisfying the above identities. We may say that a higher Lie algebroid in the sense of Definition~\ref{def.nonlinalg} gives a special example of these `two-layer' objects where the `first-layer Lie algebra' is the nilpotent Lie superalgebra of vector fields of negative weights on $F$, which is entirely fixed (locally) by the graded dimension of $F$. Is then a `second layer' (abstractly defined) always of the form introduced in this section, i.e., specified by some homological field of weight $+1$\,? (Our conjecture is, yes.  We hope to explore that elsewhere.) It may be interesting to have a look at a simple example.

\section{A model example}

We shall slightly change the notation. Let $y^i$ be variables of weight $1$ and $z^{\mu}$ be variables of weight $2$. We do not consider variables of weight $0$, so what follows should be regarded as the simplest case of a ``non-linear Lie superalgebra'' (a non-linear Lie algebroid over a point). Let $F$ be a graded manifold with coordinates $y^i,z^{\mu}$. We allow arbitrary polynomial transformations preserving weight, so we   have

\begin{equation}\label{eq.coordchange}
\begin{aligned}
    y^i&=y^{i'}T_{i'}^i\,, \\
    z^{\mu}&=z^{\mu'}T_{\mu'}^{\mu}+\frac{1}{2}\,y^{i'}y^{j'}T_{j'i'}^{\mu}\,.
\end{aligned}
\end{equation}
Note that among $y^i,z^{\mu}$ there may be  even and odd variables. The space $\Vectn(F)=\Vect_{-2}(F)\oplus \Vect_{-1}(F)$ is spanned by the vector fields
\begin{equation*}
    e_i=\der{}{y^i}\,, \quad   e_{\mu}^i=y^i\der{}{z^{\mu}}\,, \quad \text{and}\quad  e_{\mu}=\der{}{z^{\mu}}\,,
\end{equation*}
where $e_i$, $e_{\mu}^j$ make a basis of $\Vect_{-1}(F)$ and $e_{\mu}$ make a basis of $\Vect_{-2}(F)$. They satisfy
\begin{equation}\label{eq.commut}
\begin{aligned}%
 \vphantom{[}[e_i,e_j] &=0\,, \quad [e_i,e^j_{\mu}]=\delta_i^je_{\mu}\,, \quad [e^i_{\mu},e^j_{\nu}]=0\,,\\
    [e_i,e_{\mu}]&=0\,, \quad [e^i_{\mu},e_{\nu}]=0\,,\\
    [e_{\mu},e_{\nu}]&=0\,.
\end{aligned}
\end{equation}
%\end{equation*}
Under a coordinate change~\eqref{eq.coordchange}, the basis elements $e_i$, $e_{\mu}^j$, and $e_{\nu}$ transform as follows:
\begin{gather}
e_{i'}=T_{i'}^i\,e_i+(-1)^{\jt(\itt'+\mut+1)}T_j^{j'}T_{j'i'}^{\mu}\,e_{\mu}^j\,, \quad
e_{\mu'}^{i'}=(-1)^{\itt(\itt +\itt'+\mut+\mut')} T_i^{i'}T_{\mu'}^{\mu}\,e^i_{\mu}\,, \label{eq.translaw1}\\
e_{\mu'}=T_{\mu'}^{\mu}\,e_{\mu}\label{eq.translaw2}
\end{gather}
(signs are not really important and are shown only for completeness).

A vector field $Q$ of weight $1$ on $F$ has the general form
\begin{multline*}
    Q=Q^k\der{}{y^k}+Q^{\la}\der{}{z^{\la}}=\left(z^{\mu}Q_{\mu}^k+\frac{1}{2}y^iy^jQ_{ji}^k\right)\der{}{y^k}
    +\left(y^iz^{\mu}Q_{\mu i}^{\la}+\frac{1}{3!}y^iy^jy^lQ_{lji}^{\la}\right)\der{}{z^{\la}}\,.
\end{multline*}
If we ignore weights, it contains terms of degrees $1$, $2$, and $3$ in the variables $y^i,z^{\mu}$. Suppose $Q$ is odd and satisfies $Q^2=0$. Then it defines on the Lie superalgebra $\Vectn(F)$, with a basis   $e_i$, $e_{\mu}^j$, $e_{\nu}$, a `second layer', i.e., operations $d$ and $\{\_,\_\}$, given by the formulas~\eqref{eq.defd1}, \eqref{eq.defbracket1}. By a direct calculation we obtain (in the subsequent formulas we suppress exact signs because they are not relevant for our current purpose)
\begin{equation}
    de_i=0\,, \quad de^i_{\mu}=0\,, \quad de_{\mu}=\pm Q_{\mu}^ke_k\, \pm Q_{\mu i}^{\la}e_{\la}^i\,,
\end{equation}
and
\begin{equation}\label{eq.newbrack}
\begin{aligned}
    \{e_i,e_j\}&=\pm Q_{ij}^ke_k \,\pm Q_{ijk}^{\la}e^k_{\la}\,,\\
    \{e_i,e^j_{\mu}\}&=\pm Q_{ik}^je^k_{\mu}\,\pm Q_{\mu i}^{\la}e_{\la}^j\,,\\
    \{e_i,e_{\mu}\}&=\pm Q_{\mu i}^{\la}e_{\la}\,,\\
    \{e_{\mu},e_{\nu}\}&=0\,.
\end{aligned}
\end{equation}
The calculation of the remaining brackets: $\{e^i_{\mu},e_j\}$, $\{e^i_{\mu},e^j_{\nu}\}$, $\{e^i_{\mu},e_{\nu}\}$, $\{e_{\mu}, e_{i}\}$, and $\{e_{\mu},e^j_{\nu}\}$, is left to the reader. We see that the coefficients of the expansion of the homological  field $Q$ appear in various combinations as structure constants of the operations  $d$ and $\{\_,\_\}$. Claim: the vector field $Q$ can be recovered completely from the knowledge of these operations together with the Lie bracket $[\_,\_]$, i.e., from the `two-layer Lie algebra' structure on the graded vector space $\Vectn(F)$.

To see that, consider the quotient space of $\Vectn(F)$ w.r.t. the subspace spanned by $e^i_{\mu}$. Its elements, i.e., the equivalence classes of elements of $\Vectn(F)$, will be denoted by bar so that $\bar u$ stands for the class of $u$. A basis in this quotient is given by $\bar e_i, \bar e_{\mu}$. From the transformation laws~\eqref{eq.translaw1},\eqref{eq.translaw2} we see that the quotient is well-defined. In fact, it can be identified with the tangent space to $F$ at the origin. (In the general setup when variables of weight $0$ are present, we obtain the normal bundle to the zero section   $M\hookrightarrow F$.) Choose some lifting, i.e., a section of the projection. In fixed coordinates or at least allowing only linear transformations of variables, and ignoring   weights, we may use  the homological field $Q$ to define an $L_{\infty}$-algebra structure   on ``vector fields with constant coefficients'' (linear combinations of $e_i,e_{\mu}$) by higher derived bracket formulas (see~\cite{tv:higherder}):
\begin{equation*}
    (u_1,\ldots,u_k)=[\ldots[Q,u_1],\ldots, u_k](0)\,.
\end{equation*}
Conversely, the field $Q$ is completely defined by this $L_{\infty}$-structure, because its structure constants are exactly the coefficients of the Taylor expansion of $Q$ at the origin.  The subspace of ``vector fields with constant coefficients'' is not invariant under non-linear changes of coordinates. Better speak about the above quotient, which is invariant, but to define the $L_{\infty}$ operations a choice of lifting is needed. On the quotient space,
\begin{equation*}
    (\bar u_1,\ldots,\bar u_k)= \overline{[\ldots[Q,u_1],\ldots, u_k]}  \,,
\end{equation*}
where, for  $\bar u$ in the quotient, $u$ stands for its lift  to $\Vectn(F)$. Now, a simple but crucial observation is that all the homotopy Lie brackets $(\_, \ldots,\_)$ are expressed in terms of $d$, $\{\_,\_\}$,  and $[\_,\_]$\,. In our example we have only three operations: the differential $\p=(\_)$, the binary bracket $(\_,\_)$ and the ternary bracket $(\_,\_,\_)$, because the field $Q$ does not contain terms higher than   cubic. We have
\begin{equation*}
    \p \bar u = \overline{du}\,,\quad  (\bar u,\bar v)=\overline{\{u,v\}}\,, \quad (\bar u,\bar v,\bar w)=\overline{[\{u,v\},w]}\,.
\end{equation*}
In general, we would have $(\bar u_1,  \ldots,\bar u_k)=\overline{[\ldots [\{u_1,u_2\},u_3], \ldots, u_k]}$. Since, in fixed coordinates, this (non-invariant) $L_{\infty}$-algebra structure completely defines $Q$ and is in turn completely defined by $d$, $\{\_,\_\}$,  and $[\_,\_]$, we see that, indeed, the homological vector field $Q$ is completely determined by the structure it induces on $\Vectn(F)$.

For an illustration, in our example we have, on the basis vectors $\bar e_i,\bar e_{\mu}$,
\begin{gather*}
    \p \bar e_i=\overline{de_i}=0\,, \\
    \p \bar e_{\mu}=\overline{d e_{\mu}}=\pm Q_{\mu}^k\bar e_k\,,\\
    (\bar e_i,\bar e_j)=\overline{\{e_i,e_j\}}=\pm Q_{ij}^k\bar e_k\,,\\
    (\bar e_i,\bar e_{\mu})=\overline{\{e_i,e_{\mu}\}}=\pm Q_{\mu i}^{\la}\bar e_{\la}\,,\\
    (\bar e_{\mu},\bar e_{\nu})=\overline{\{e_{\mu},e_{\nu}\}}=0\,,\\
    (\bar e_i,\bar e_j, \bar e_k)=\overline{[\{e_i,e_j\},e_k]}=\pm Q_{ijk}^{\la}\bar e_{\la}\,,\\
    (\bar e_i,\bar e_j, \bar e_{\mu})=\overline{[\{e_i,e_j\},e_{\mu}]}=0\,,\\
    (\bar e_i,\bar e_{\mu},\bar e_{\nu})=\overline{[\{e_i,e_{\mu}\},e_{\nu}]}=0\,,\\
    (\bar e_{\mu},\bar e_{\nu},e_{\la})=\overline{[\{e_{\mu},e_{\nu}\},e_{\la}]}=0\,,
\end{gather*}
which follow by combining ~\eqref{eq.newbrack} and \eqref{eq.commut}. So all the coefficients of the expansion of $Q$ are recovered, as claimed.

\section{Conclusions}

We have  studied an algebraic structure corresponding to a ``non-linear'' or ``higher'' Lie algebroid defined geometrically as a non-negatively graded $Q$-manifold with the homological vector field of weight $+1$. For $Q$-manifolds corresponding to ordinary Lie algebroids, this would correspond to the structure on sections. In the non-linear case we use vector fields of negative weight as a replacement of Lie algebroid sections. We have established algebraic identities satisfied by new operations that arise. It is possible to treat the problem in a more abstract way in the framework of graded Lie superalgebras and the resulting new algebraic structure allows an axiomatization. From the geometric viewpoint, to the initial graded manifold (which is a non-linear fiber bundle) corresponds an associated graded vector bundle whose sections carry the described structure (which we call a `two-layer Lie algebra' or  `two-layer Lie algebroid' structure because it basically consists of a differential and  a new binary bracket sitting on  top of a Lie superalgebra). This should be explored further. We have shown, using a  model example, that in our setting hides also an $L_{\infty}$-structure, which is, however, defined non-invariantly. Still, it is in some sense equivalent to the invariant `two-layer Lie' structure, which allows us to conclude that the original geometric structure of a non-negatively graded $Q$-manifold (a ``non-linear Lie algebroid'') and the resulting new algebraic structure on vector fields (which may be defined by axioms) are in one-to-one correspondence. We hope to elaborate this elsewhere.

%%%%%%%%%%%%%%%%%%%%%%%%%%%%%%%%%%%%%%%%%%%%%%%%
%% BACKMATTER
%%%%%%%%%%%%%%%%%%%%%%%%%%%%%%%%%%%%%%%%%%%%%%%%

\begin{theacknowledgments}
I am deeply grateful to Dr Kirill Mackenzie who read this text and made many important comments, and for many discussions over years.
  It is a pleasure to thank the organizers of the annual international Workshops on
  Geometric Methods in Physics in Bia{\l}owie\.{z}a where part of this work
  was reported. I am particularly grateful to Professor Anatol~Odzijewicz for the
  hospitality and the exceptionally inspiring atmosphere at the
  meetings over so many years.
\end{theacknowledgments}

%%%%%%%%%%%%%%%%%%%%%%%%%%%%%%%%%%%%%%%%%%%%%%%%
%% The bibliography can be prepared using the BibTeX program or
%% manually.
%%
%% The code below assumes that BibTeX is used.  If the bibliography is
%% produced without BibTeX comment out the following lines and see the
%% aipguide.pdf for further information.
%%
%% For your convenience a manually coded example is appended
%% after the \end{document}
%%%%%%%%%%%%%%%%%%%%%%%%%%%%%%%%%%%%%%%%%%%%%%%%

%%%%%%%%%%%%%%%%%%%%%%%%%%%%%%%%%%%%%%%%%%%%%%%%
%% You may have to change the BibTeX style below, depending on your
%% setup or preferences.
%%
%%
%% For The AIP proceedings layouts use either
%%%%%%%%%%%%%%%%%%%%%%%%%%%%%%%%%%%%%%%%%%%%
\bibliographystyle{plain}
%\bibliographystyle{aipproc}   % if natbib is available
%\bibliographystyle{aipprocl} % if natbib is missing

%%%%%%%%%%%%%%%%%%%%%%%%%%%%%%%%%%%%%%%%%%%
%% You probably want to use your own bibtex database here
%%%%%%%%%%%%%%%%%%%%%%%%%%%%%%%%%%%%%%%%%%%

%\bibliography{geometry}
%\end{document}

\def\cprime{$'$}

\end{document}